# Uniform limit laws of the logarithm for estimators of the additive regression function in the presence of right censored data

## Mohammed Debbarh

*Laboratoire de Statistique Thérique et Appliquée, Université Paris VI, 175 rue du Chevaleret 75013 Paris France*
*e-mail:* debbarh@ccr.jussieu.fr

and

## Vivian Viallon[*]

*Department of Biostatistics, Hôpital Cochin, APHP, Faculté de Médecine, Université Paris Descartes, Paris, France*
*e-mail:* vivian.viallon@univ-paris5.fr

**Abstract:** It has been recently shown that nonparametric estimators of the additive regression function could be obtained in the presence of right censoring by coupling the marginal integration method with initial kernel-type Inverse Probability of Censoring Weighted estimators of the multivariate regression function [10]. In this paper, we get the exact rate of strong uniform consistency for such estimators. Our uniform limit laws especially lead to the construction of asymptotic simultaneous 100% confidence bands for the true regression function.

## Contents

*Correspondence to: Vivian Viallon Dⁱpartement de Biostatistique, Hôpital Cochin, 27 rue du Faubourg Saint Jacques, 75014 Paris, France. Phone: + 33 1 58 41 19 60. Fax: +33 1 58 41 19 61. E-mail: vivian.viallon@univ-paris5.fr

## 1. Introduction

Consider a triple $(Y, C, \mathbf{X})$ of random variables defined in $\mathbb{R}^+ \times \mathbb{R}^+ \times \mathbb{R}^d$, $d \geq 2$, where $Y$ is the variable of interest (typically a lifetime variable), $C$ a censoring variable and $\mathbf{X} = (X_1, \ldots, X_d)$ a vector of concomitant variables. In most practical applications, such as epidemiology or reliability, the relationship between $Y$ and $\mathbf{X}$ is of particular interest. Denoting by $\psi$ a given measurable function, we will focus here on the study of the conditional expectation of $\psi(Y)$ given $\mathbf{X} = \mathbf{x}$,

$$m_\psi(\mathbf{x}) = \mathbb{E}\big(\psi(Y) \mid \mathbf{X} = \mathbf{x}\big), \quad \text{for all } \mathbf{x} = (x_1, \ldots, x_d) \in \mathbb{R}^d. \tag{1.1}$$

The introduction of the function $\psi$ will allow us to treat simultaneously the standard regression function and the conditional distribution function (see Remark 2.4 below).

In the right censorship model, the pair $(Y, C)$ is not directly observed and the corresponding information is given by $Z = \min\{Y, C\}$ and $\delta = \mathbb{1}_{\{Y \leq C\}}$, $\mathbb{1}_E$ standing for the indicator function of the set $E$. Therefore, we will assume that a sample $\mathcal{D}_n = \{(Z_i, \delta_i, \mathbf{X}_i), i = 1, \ldots, n\}$ of independent and identically distributed replicae of the triple $(Z, \delta, \mathbf{X})$ is at our disposal. In this setting, transformations of the observed data $\mathcal{D}_n$ are usually needed to estimate functionals of the conditional law of $Y$ (see, e.g., [4, 17, 28–30] and the recent work of [35]). Estimators based on these transformations are usually referred to as *synthetic data* estimators. In this paper, following the ideas initiated by [28], we use a nonparametric version of particular synthetic data estimators, commonly referred to as *Inverse Probability of Censoring Weighted* [*I.P.C.W.*] estimators (see [3, 5] and [27] for some results related to nonparametric I.P.C.W. estimates of the censored regression function). It is however noteworthy that the methodology we propose here for *I.P.C.W.*-type estimators shall apply with minor modifications to cope with other synthetic data estimators (see Paragraph 3.1 below).

A well-known issue in nonparametric estimation is the so-called *curse of dimensionality*: the rate of convergence of nonparametric estimators generally decreases as the dimensionality $d$ of the covariate increases. To get round this problem, one solution is to work, if possible, under the additive model assumption, which allows to write the regression function as follows,

$$m_\psi(\mathbf{x}) = \mu + \sum_{\ell=1}^{d} m_{\psi,\ell}(x_\ell). \tag{1.2}$$

In (1.2), the real-valued functions $m_{\psi,\ell}$, $\ell = 1, \ldots, d$, are defined up to an additive constant, and the assumption $\mathbb{E} m_{\psi,\ell}(X_\ell) = 0$, $\ell = 1, \ldots, d$, is usually made to ensure identifiability. This assumption implies $\mu = \mathbb{E}\psi(Y)$. In the uncensored case, several methods have been proposed to estimate the additive regression function. We shall evoke, among others, the methods based on *B*-splines [36], on the *backfitting* algorithm [21, 23, 32] and on marginal integration [34, 40, 31]. In [17], Fan and Gijbels established the asymptotic normality for estimators obtained *via* the backfitting algorithm combined with various synthetic data estimators. In [3], Brunel and Comte considered additive models as special cases in their study of adaptive projection I.P.C.W. estimators. Here, following the ideas introduced in [10], we make use of the marginal integration method, coupled with initial kernel-type I.P.C.W. estimators to provide an estimator for the additive censored regression function. This combination leads to estimators for which the theory is easier to derive, which was wanted here, given the technicalities in the proof, even in this simplified setting (note however that, as already mentioned, extensions to other synthetic data estimators can be obtained; see Paragraph 3.1). In a previous work [10], the mean-square convergence rate was established for the integrated estimator defined in (2.7) below. In the present paper, we get the exact corresponding rate of strong uniform consistency (see Theorem 3.2 below). Our limit law corresponds to the extension of Theorem 2 in [9] to the censored case. Moreover, following the ideas developed in [13], asymptotic simultaneous 100% confidence bands are derived for the true regression function. This kind of bands may be complementary to the more classical $(1 - \alpha) \times 100\%$ pointwise confidence intervals derived from CLT type results (see Section 4).

## 2. Hypotheses-Notations

Before presenting our estimator and stating our results, we shall introduce some notations as well as our working assumptions. First consider the hypotheses to be made on the random triple $(Y, C, \mathbf{X})$. Introduce, for all $t \in \mathbb{R}$, $F(t) = \mathbb{P}(Y > t)$, $G(t) = \mathbb{P}(C > t)$ and $H(t) = \mathbb{P}(Z > t)$, the right continuous survival functions pertaining to $Y$, $C$ and $Z$ respectively.

$(C.1)$   $C$ and $Y$ are independent and $\mathbb{P}(Y \leq C | \mathbf{X}, Y) = \mathbb{P}(Y \leq C | Y)$.
$(C.2)$   $G$ is continuous.
$(C.3)$   $m_\psi$ is $s$-times continuously differentiable, $s \geq 1$, and

$$\sup_{\mathbf{x}} \left| \frac{\partial^s}{\partial x_1^{s_1} \ldots \partial x_d^{s_d}} \, m_\psi(\mathbf{x}) \right| < \infty; \; s_1 + \cdots + s_d = s.$$

**Remark 2.1.** *Assumption $(C.3)$ will allow to control bias terms. Assumptions $(C.1)$ and $(C.2)$ are essentially needed when using most synthetic data estimators. $(C.2)$ allows to use convergence results for the Kaplan-Meier [25] estimator of $G$. In addition, $(C.1)$ especially allows to derive the result (2.1) below, which is a fundamental requirement for synthetic data. This assumption was also used by Stute [38] in another context. It is however noteworthy that Beran [2] (see* *also [8] and [12]) worked under the weaker assumption of conditional independence between $Y$ and $C$ given $\mathbf{X}$ to derive properties for a local version of the Kaplan-Meier estimator. On the other hand, to use Beran's local Kaplan-Meier estimator the censoring has to be locally fair, i.e., such that $\mathbb{P}(C > t|\mathbf{X}) > 0$ whenever $\mathbb{P}(Y > t|\mathbf{X})$. Here, (see assumption $(\mathbf{A})(ii)$ below), we essentially suppose that $G(t) > 0$ whenever $F(t) > 0$, which is on its turn a weaker assumption. For a nice discussion on the differences between the assumptions to be made when using either Beran's estimator or I.P.C.W. type estimators, we refer to [5].*

Denote by $f$ [resp. $f_\ell$, $\ell = 1, \ldots, d$] the density of $\mathbf{X}$ [resp. $X_\ell$, $\ell = 1, \ldots, d$]. Further let $\mathcal{C}_1, \ldots, \mathcal{C}_d$, be $d$ compact intervals of $\mathbb{R}$ with non empty interior, and set $\mathcal{C} = \mathcal{C}_1 \times \cdots \times \mathcal{C}_d$ the corresponding product. For every subset $\mathcal{E}$ of $\mathbb{R}^q$, $q \geq 1$, and any $\alpha > 0$, introduce the $\alpha$-neighborhood $\mathcal{E}^\alpha$ of $\mathcal{E}$,

$$\mathcal{E}^\alpha = \{x : \inf_{y \in \mathcal{E}} |x - y|_{\mathbb{R}^q} \leq \alpha\},$$

with $|\cdot|_{\mathbb{R}^q}$ standing for the usual euclidian norm on $\mathbb{R}^q$. The functions $f$ and $f_\ell$, $\ell = 1, \ldots, d$, will be supposed to be continuous, and we will assume the existence of a constant $\alpha > 0$ such that the following assumptions hold.

($C$.4)  $\forall x_\ell \in \mathcal{C}_\ell^\alpha, f_\ell(x_\ell) > 0$, $\ell = 1, \ldots, d$, and $\forall \mathbf{x} \in \mathcal{C}^\alpha, f(\mathbf{x}) > 0$.
($C$.5)  $f$ is $s'$-times continuously differentiable on $\mathcal{C}^\alpha, s' > sd$.

**Remark 2.2.** *Assumption $(C$.4) is classical when dealing with kernel type estimators of the regression function (see, e.g., [13, 15]). The fact that $s' > sd$ in $(C$.5), when combined with $(C$.3) above and $(K$.1-2) and $(H$.4) below, allows to derive easily the results pertaining to the case where the density function $f$ is unknown from the ones obtained in the simpler case where this function is known. Some refinements in our proofs might allow for relaxing $(C$.5) (see, e.g., [22]).*

Recalling (1.1), we will let $\psi$ vary in a *pointwise measurable VC subgraph class* $\mathcal{F}$ of measurable real-valued functions defined on $\mathbb{R}$ (for the definitions of pointwise measurable classes of functions and VC subgraph classes of functions, we refer to p. 110 and Chapter 2 in [41]). We will also assume that $\mathcal{F}$ has a measurable envelope function $\Upsilon(y) \geq \sup_{\psi \in \mathcal{F}} |\psi(y)|$, $y \in \mathbb{R}$, such that

($C$.6)  $\Upsilon$ is uniformly bounded on $\mathbb{R}$.

**Remark 2.3.** *In the uncensored setting, $(C$.6) can be replaced by some finiteness condition on the moment of order 2 of $\Upsilon(Y)$ (see [13] or [15]). In the censored setting however, such refinements are useless due to the assumption $(\mathbf{A})$ below.*

**Remark 2.4.** *Choices of particular interest for the class $\mathcal{F}$ are $\mathcal{F}_{reg} = \{I\}$, where $I$ denotes the identity function on $\mathbb{R}$ and $\mathcal{F}_{cdf} = \{\mathbb{1}_{(-\infty, t]}, t \in \mathbb{R}\}$. Considering the class $\mathcal{F}_{reg}$ allows to treat the case of the classical regression function. On the other hand, considering the class $\mathcal{F}_{cdf}$ allows to derive the uniform consistency (especially over $t \in \mathbb{R}$) for estimates of the conditional distribution function. We refer to [15] for examples in the uncensored case.*

We will further employ sequences of positive constants $\{h_n\}_{n \geq 1}$ and $\{h_{\ell,n}\}_{n \geq 1}$, $1 \leq \ell \leq d$, satisfying the following conditions.

$(H.1)$  $h_n \downarrow 0$, $h_{\ell,n} \downarrow 0$, $nh_n^d \uparrow \infty$ and $nh_{\ell,n} \uparrow \infty$ as $n \to \infty$.

$(H.2)$  $nh_n^d / \log n \to \infty$ and $nh_{\ell,n} / \log n \to \infty$ as $n \to \infty$.

$(H.3)$  $nh_{\ell,n} \prod_{j=1}^d h_{j,n}^{2s_j} / |\log h_{\ell,n}| \to 0$, for all $s_1 + \cdots + s_d = s$, as $n \to \infty$.

$(H.4)$  $h_{\ell,n} \log n / (h_n^d |\log h_{\ell,n}|) \to 0$ as $n \to \infty$.

$(H.5)$  $\log \log n / |\log h_n| \to 0$ and $\log \log n / |\log h_{\ell,n}| \to 0$ as $n \to \infty$.

**Remark 2.5.** *Assumptions $(H.1\text{-}2\text{-}5)$ are classical in the empirical process theory, and are often referred to as the Csörgő-Révész-Stute [CRS] conditions [7, 37]. They especially allow to control variance-type terms. On the other hand, assumption $(H.3)$ allows to control bias terms (see Lemma 5.8 below). As already mentioned, assumption $(H.4)$ allows to derive easily the results pertaining to the case where the density function $f$ is unknown from the ones obtained in the simpler case where this function is known.*

As mentioned in [19], functionals of the (conditional) law can generally not be estimated on the complete support when the variable of interest is right-censored. So, to state our results, we will work under the assumption $(\mathbf{A})$, that will be said to hold if either $(\mathbf{A})(i)$ or $(\mathbf{A})(ii)$ below holds. For any right continuous survival function $L$ defined on $\mathbb{R}$, set $T_L = \sup\{t \in \mathbb{R} : L(t) > 0\}$.

$(\mathbf{A})(i)$  There exists a $\omega < T_H$ such that, for all $\psi \in \mathcal{F}$, $\psi = 0$ on $(\omega, \infty)$.

$(\mathbf{A})(ii)$  $(a)$  For a given $0 < p \leq 1/2$, $\int_0^{T_H} -F^{-p/(1-p)} dG < \infty$;
$\qquad\quad\ (b)$  $T_F < T_G$;
$\qquad\quad\ (c)$  $n^{2p-1} h_{\ell,n}^{-1} |\log(h_{\ell,n})| \to \infty$, as $n \to \infty$, for every $\ell = 1, \ldots, d$.

It is noteworthy that the assumption $(\mathbf{A})(ii)$ is needed in our proofs when considering the estimation of the "classical" regression function, which corresponds to the choice $\psi(y) = y$. On the other hand, rates of convergence for estimators of functionals such as the conditional distribution function $\mathbb{P}(Y \leq t | \mathbf{X})$ can be obtained under weaker conditions, when restricting ourselves to $t \in [0, \omega]$ with $\omega < T_H$.

These preliminaries being given, we can recall the procedure we proposed in [10] to estimate the censored regression function under the additive model assumption. Let $K$ be a bounded and compactly supported kernel on $\mathbb{R}^d$. By kernel, we mean as usual a measurable function integrating to one on its support. We define the kernel density estimator $\hat{f}_n$ of $f$ by

$$\hat{f}_n(\mathbf{x}) = \frac{1}{nh_n^d} \sum_{i=1}^d K\left(\frac{\mathbf{x} - \mathbf{X}_i}{h_n}\right).$$

Now, as was observed notably by Koul et al. [28], we have under $(C.1)$,

$$\mathbb{E}\left\{\frac{\delta\psi(Z)}{G(Z)} \mid \mathbf{X}\right\} = \mathbb{E}\left\{\frac{\psi(Y)}{G(Y)}\mathbb{E}\left(\mathbb{1}_{\{Y \leq C\}} | \mathbf{X}, Y\right) \mid \mathbf{X}\right\} = \mathbb{E}(\psi(Y) | \mathbf{X}). \quad (2.1)$$

Then, denoting by $G_n^\star$ the Kaplan-Meier [25] estimator of $G$, kernel-type estimators of the multivariate regression function $m_\psi(\mathbf{x})$ defined in (1.1) can be easily constructed [27]. Here, because marginal integration will further be applied, the internal estimator idea of Jones [24] has to be used. That leads us to consider the following multivariate I.P.C.W. kernel-type estimator of the regression function,

$$\widetilde{m}_{\psi,n}^\star(\mathbf{x}) = \frac{1}{n} \sum_{i=1}^n \left\{ \frac{\delta_i \psi(Z_i)}{G_n^\star(Z_i)\hat{f}_n(\mathbf{X}_i)} \prod_{\ell=1}^d \frac{1}{h_{\ell,n}} K_\ell\left(\frac{x_\ell - X_{i,\ell}}{h_{\ell,n}}\right) \right\}. \tag{2.2}$$

Here the kernel functions $K_\ell$, $\ell = 1, \ldots, d$, defined in $\mathbb{R}$ are supposed to be continuous, of bounded variation (i.e. such that $0 < \int_{\mathbb{R}} |dK_\ell(t)| < \infty$) and compactly supported. Recalling that a kernel function $\Gamma$ defined in $\mathbb{R}^d$ is said to be of order $\gamma$, for any $\alpha \geq 1$, whenever (*a*) and (*b*) below holds jointly,

(*a*) $\int_{\mathbb{R}^d} u_1^{j_1} \ldots u_d^{j_d} \Gamma(\mathbf{u}) d\mathbf{u} = 0$, $j_1, \ldots, j_d \geq 0$, $j_1 + \cdots + j_d = 0, 1, \ldots, \gamma - 1$;
(*b*) $\int_{\mathbb{R}^d} |u_1^{j_1} \ldots u_d^{j_d}| |\Gamma(\mathbf{u})| d\mathbf{u} < \infty$, $j_1, \ldots, j_d \geq 0$, $j_1 + \cdots + j_d = \gamma$;

we will also impose the conditions (*K*.1-2).

(*K*.1) $\mathbb{K} := \prod_{\ell=1}^d K_\ell$ is of order $s$.
(*K*.2) $K$ is of order $s'$.

In order to apply the marginal integration method (see [31, 34]), introduce $q_1, \ldots, q_d$, $d$ given density functions defined in $\mathbb{R}$. Further set, for all $\mathbf{x} = (x_1, .., x_d) \in \mathbb{R}^d$, $q(\mathbf{x}) = \prod_{\ell=1}^d q_\ell(x_\ell)$ and, for every $\ell = 1, \ldots, d$, $q_{-\ell}(\mathbf{x}_{-\ell}) = \prod_{j \neq \ell} q_j(x_j)$ with $\mathbf{x}_{-\ell} = (x_1, .., x_{\ell-1}, x_{\ell+1}, .., x_d)$. Now, we can define

$$\eta_{\psi,\ell}(x_\ell) = \int_{\mathbb{R}^{d-1}} m_\psi(\mathbf{x}) q_{-\ell}(\mathbf{x}_{-\ell}) d\mathbf{x}_{-\ell} - \int_{\mathbb{R}^d} m_\psi(\mathbf{x}) q(\mathbf{x}) d\mathbf{x}, \quad \ell = 1, \ldots, d, \tag{2.3}$$

in such a way that, recalling (1.2), the two following equalities hold,

$$\eta_{\psi,\ell}(x_\ell) = m_{\psi,\ell}(x_\ell) - \int_{\mathbb{R}} m_{\psi,\ell}(u) q_\ell(u) du, \quad \ell = 1, \ldots, d, \tag{2.4}$$

$$m_\psi(\mathbf{x}) = \sum_{\ell=1}^d \eta_{\psi,\ell}(x_\ell) + \int_{\mathbb{R}^d} m_\psi(\mathbf{u}) q(\mathbf{u}) d\mathbf{u}. \tag{2.5}$$

In view of (2.4) and (2.5), for every $\ell = 1, \ldots, d$, $\eta_{\psi,\ell}$ and $m_{\psi,\ell}$ are equal up to an additive constant, so that the functions $\eta_{\psi,\ell}$ are actually some additive components, which coincide with $m_{\psi,\ell}$ for the choice $q_\ell = f_\ell$ (which is only achievable if $f_\ell$ is known). From (2.2) and (2.3), a natural estimator of the $\ell$-th component $\eta_{\psi,\ell}$ is given, for $\ell = 1, \ldots, d$, by

$$\hat{\eta}_{\psi,\ell}^\star(x_\ell) = \int_{\mathbb{R}^{d-1}} \widetilde{m}_{\psi,n}^\star(\mathbf{x}) q_{-\ell}(\mathbf{x}_{-\ell}) d\mathbf{x}_{-\ell} - \int_{\mathbb{R}^d} \widetilde{m}_{\psi,n}^\star(\mathbf{x}) q(\mathbf{x}) d\mathbf{x}. \tag{2.6}$$

From (2.5) and (2.6), an estimator $\widehat{m}_{\psi,add}^{\star}$ of the additive regression function can be deduced,

$$\widehat{m}_{\psi,add}^{\star}(\mathbf{x}) = \sum_{\ell=1}^{d} \widehat{\eta}_{\psi,\ell}^{\star}(x_\ell) + \int_{\mathbf{R}^d} \widetilde{m}_{\psi,n}^{\star}(\mathbf{u}) q(\mathbf{u}) d\mathbf{u}. \qquad (2.7)$$

In the sequel, we will assume that the known integration density function $q_\ell$ has a compact support included in $\mathcal{C}_\ell$, $\ell = 1, \ldots, d$. Moreover, we will impose the following assumption on the functions $q_{-\ell}$, $\ell = 1, \ldots, d$.

$(Q.1)$     $q_{-\ell}$ is a bounded and $s$-times differentiable function such that

$$\sup_{\mathbf{x}_{-\ell}} \left| \frac{\partial^s}{\partial^{i_1} x_1 \ldots \partial^{i_d} x_d} q_{-\ell}(\mathbf{x}_{-\ell}) \right| < \infty, \ i_1 + \cdots + i_d = s, \ \ell = 1, \ldots, d.$$

Before stating our main results, some additional notations are needed. For all $\psi \in \mathcal{F}$, all $\mathbf{u} = (u_1, \ldots, u_d) \in \mathcal{C}$ and every $\ell = 1, \ldots, d$, set

$$H_\psi(\mathbf{u}) = \mathbb{E}\left( \frac{\psi^2(Y)}{G(Y)} \Big| \mathbf{X} = \mathbf{u} \right) \qquad (2.8)$$

$$\text{and} \quad \phi_{\psi,\ell}(u_\ell) = \int_{\mathbf{R}^{d-1}} \frac{H_\psi(\mathbf{u})}{f(\mathbf{u}_{-\ell}|u_\ell)} q_{-\ell}(\mathbf{u}_{-\ell}) d\mathbf{u}_{-\ell}. \qquad (2.9)$$

Further set, for all $\psi \in \mathcal{F}$ and every $\ell = 1, \ldots, d$,

$$\sigma_{\psi,\ell} = \sup_{x_\ell \in \mathcal{C}_\ell} \sqrt{\frac{\phi_{\psi,\ell}(x_\ell)}{f_\ell(x_\ell)} \int_{\mathbf{R}} K_\ell^2}, \quad \sigma_\ell = \sup_{\psi \in \mathcal{F}} \sigma_{\psi,\ell} \ \text{ and } \ \sigma = \sum_{\ell=1}^{d} \sigma_\ell. \qquad (2.10)$$

## 3. Main results

We have now all the ingredients to state our results. From now on, $\xrightarrow{a.s.}$ will stand for almost sure convergence. Theorem 3.1 below describes the asymptotic behavior of the additive component estimates $\widehat{\eta}_{\psi,\ell}$, $\ell = 1, \ldots, d$, defined in (2.6).

**Theorem 3.1.** *Under the hypotheses* (**A**), $(C.1\text{-}2\text{-}3\text{-}4\text{-}5\text{-}6)$, $(H.1\text{-}2\text{-}3\text{-}4\text{-}5)$, $(K.1\text{-}2)$ *and* $(Q.1)$, *we have, for* $\ell = 1, \ldots, d$,

$$\sqrt{\frac{nh_{\ell,n}}{2|\log h_{\ell,n}|}} \sup_{\psi \in \mathcal{F}} \sup_{x_\ell \in \mathcal{C}_\ell} \pm \{\widehat{\eta}_{\psi,\ell}^{\star}(x_\ell) - \eta_{\psi,\ell}(x_\ell)\} \xrightarrow{a.s.} \sigma_\ell \quad \text{as } n \to \infty, \qquad (3.1)$$

*where* $\sigma_\ell$ *is as in* (2.10).

From Theorem 3.1, we will deduce an analogous result for the additive regression function estimator $\widehat{m}_{\psi,add}^{\star}$ defined in (2.7).

**Theorem 3.2.** *Assume the hypotheses of Theorem 3.1 hold. If, in addition,* $h_{\ell,n} = h_{1,n}$ *for every* $\ell = 1, \ldots, d$, *then we have,*

$$\sqrt{\frac{nh_{1,n}}{2|\log h_{1,n}|}} \sup_{\psi \in \mathcal{F}} \sup_{\mathbf{x} \in \mathcal{C}} \pm \{\widehat{m}_{\psi,add}^{\star}(\mathbf{x}) - m_\psi(\mathbf{x})\} \xrightarrow{a.s.} \sigma \quad \text{as } n \to \infty. \qquad (3.2)$$

*where* $\sigma$ *is as in* (2.10).

Keep in mind that a similar result is readily obtained for the conditional distribution function by selecting $\mathcal{F} = \{\mathbb{I}_{[0,t]}, t \in \mathbb{R}^+\}$ in Theorem 3.2.

The proofs of Theorems 3.1 and 3.2 are postponed to Section 5. A sketch of the proof of Theorem 3.1 is as follows. It will be split into two main parts. First we will assume that both the survival function $G$ of $C$ and the density function $f$ of $\mathbf{X}$ are known. Then, using appropriate approximations lemmas, we will show how to treat the general case (i.e. the case where neither $f$ nor $G$ is known). To establish the results in the case where both $G$ and $f$ are known, we will mostly borrow the arguments developed in [13] and [15] (see also [16]), which rest on recent developments in empirical process theory, and especially on an exponential bound due to Talagrand [39] (see also Inequality A.1 in the Appendix).

In the following Paragraph 3.1, we show how our results may be extended to the case of more general synthetic data. In Section 4 we present an application of our results, following the ideas developed in [13].

### *3.1. Extensions*

Here, we will limit ourselves to the case $\mathcal{F} = \{I\}$, where $I$ stand for the identity function on $\mathbb{R}$. The corresponding estimator defined in (2.2), and then the one defined under the additive assumption in (2.7), rest on the following transformation, which is due to Koul et al. [28]: for $1 \leq i \leq n$,

$$(\delta_i, Z_i) \longrightarrow \frac{\delta_i Z_i}{G_n^\star(Z_i)}, \tag{3.3}$$

which, in the case where $G$ is known, reduces to

$$(\delta_i, Z_i) \longrightarrow \frac{\delta_i Z_i}{G(Z_i)}. \tag{3.4}$$

Note that (3.4) sets a censored observation to 0 and multiplies an uncensored observation by a factor $[G(Z_i)]^{-1}$, which can be very large if $G(Z_i)$ is near 0. Alternative, and more general, synthetic data can be constructed in the following way. For any given $\rho \in \mathbb{R}$, set

$$\begin{aligned}
\Theta_1(z) &= (1+\rho) \int_0^z \frac{dt}{G(t)} - \frac{\rho z}{G(z)}, \\
\Theta_2(z) &= (1+\rho) \int_0^z \frac{dt}{G(t)},
\end{aligned} \tag{3.5}$$

with $\rho$ chosen such that $\Theta_1(Z) > 0$ almost surely. Now, consider the transformation, for $1 \leq i \leq n$,

$$(\delta_i, Z_i) \longrightarrow \widetilde{Y_i} := \delta_i \Theta_1(Z_i) + (1 - \delta_i)\Theta_2(Z_i). \tag{3.6}$$

Observe that (3.4) corresponds to the particular choice $\rho = -1$. The choice $\rho = 0$ is also popular, and was first considered in [30]. Other choices (including some data-dependent choices) are discussed in [17].

Consider the following independence assumption.

$(\widetilde{C.1})$  $C$ and $(Y, \mathbf{X})$ are independent.

Observe that $(\widetilde{C.1})$ naturally implies $(C.1)$. Under the assumption $(\widetilde{C.1})$, we have (see, e.g., [26])

$$\mathbb{E}(\widetilde{Y}_i | \mathbf{X}) = \mathbb{E}(Y_i | \mathbf{X}).$$

A close look into the proof presented in Section 5 below reveals that, in the case where $\mathcal{F} = \{I\}$, Theorems 3.1 and 3.2 still hold when considering estimators built on the "general" synthetic data, under the assumption $(\widetilde{C.1})$. The only difference is the term $H_\psi(\mathbf{u}) = H_I(\mathbf{u})$ (since $\psi = I$) defined in (2.8) that shall be replaced in the general case by

$$\widetilde{H}_I(\mathbf{u}) = \mathbb{E}(\widetilde{Y}^2 | \mathbf{X} = \mathbf{u}),$$

where $\widetilde{Y} := \delta \Theta_1(Z) + (1 - \delta)\Theta_2(Z)$.

## 4. Application

Following the ideas developed in [13], we now present a practical application of Theorem 3.1. Recall the definition (2.9) of the functions $\phi_{\psi, \ell}$. Then, for any fixed $\psi \in \mathcal{F}$, and every $\ell = 1, \ldots, d$, let $\widehat{\tau}_{\psi, \ell, n}(x_\ell)$ be a consistent estimator of $\tau_{\psi, \ell}(x_\ell)$, with $\tau_{\psi, \ell}(x_\ell) = \sqrt{\phi_{\psi, \ell}(x_\ell)/f_\ell(x_\ell)}$. For instance, set

$$
\begin{aligned}
\widehat{\tau}_{\psi, \ell, n}(x_\ell) \;=\; & \frac{1}{n h_{\ell, n}} \sum_{i=1}^n \frac{\delta_i \psi^2(Z_i)}{G_n^{*2}(Z_i)} K_\ell\left(\frac{x_\ell - X_{i, \ell}}{h_{\ell, n}}\right) \\
& \times \int_{\mathbb{R}^{d-1}} \frac{\prod_{j \neq \ell} h_{j, n}^{-1} K\left(\frac{x_j - X_{i, j}}{h_{j, n}}\right)}{\widehat{f}_n(\mathbf{x})} \; q_{-\ell}(\mathbf{x}_{-\ell}) d\mathbf{x}_{-\ell}.
\end{aligned}
$$

Further set

$$L_n(x_\ell) = \left\{\frac{2|\log h_{\ell, n}|}{n h_\ell} \times \widehat{\tau}_{\psi, \ell, n}(x_\ell)\right\}^{1/2} \left[\int_{\mathbb{R}} K_\ell^2\right]^{1/2}.$$

In view of Theorem 3.1, it is straightforward that, for each $0 < \varepsilon < 1$, there exists almost surely an $n_0 = n_0(\varepsilon)$ such that, for all $n \geq n_0$,

$$
\begin{aligned}
&\eta_{\psi, \ell}(x_\ell) \in \left[\widehat{\eta}_{\psi, \ell}^{\star}(x_\ell) \pm (1 + \varepsilon) L_n(\boldsymbol{x})\right], \quad \text{uniformly over } x_\ell \in \mathcal{C}_\ell, \\
&\eta_{\psi, \ell}(x_\ell) \notin \left[\widehat{\eta}_{\psi, \ell}^{\star}(x_\ell) \pm (1 - \varepsilon) L_n(\boldsymbol{x})\right], \quad \text{for some } x_\ell \in \mathcal{C}_\ell.
\end{aligned}
\tag{4.1}
$$

Therefore, under the assumptions of Theorem 3.1, the interval

$$\left[A_{n, \ell}(x_\ell), \quad B_{n, \ell}(x_\ell)\right] := \left[\widehat{\eta}_{\psi, \ell}^{\star}(x_\ell) - L_n(x_\ell), \quad \widehat{\eta}_{\psi, \ell}^{\star}(x_\ell) + L_n(x_\ell)\right], \tag{4.2}$$

provides *asymptotic simultaneous confidence bands* (at an asymptotic confidence level of 100 %) for $\eta_{\psi, \ell}(x_\ell)$ over $x_\ell \in \mathcal{C}_\ell$ (see [13] for more details). It is noteworthy that our bands do not provide confidence regions in the usual sense, since they are not based on a specified confidence level $1 - \alpha$. Instead, they hold with probability tending to 1 as $n \to \infty$, and are then more conservative (since they are simultaneous and with an asymptotic level of 100 %). A comparison between pointwise $(1 - \alpha) \times 100\%$ confidence intervals and our simultaneous almost certainty bands can be found in [11]. In most applications, we recommend the construction of both types of confidence region to assess the form of the relationship between $\psi(Y)$ and $\mathbf{X}$.

**Remark 4.1.** *For finite sample size use, Deheuvels and Mason [13] give some recommendations on how to ensure that the simultaneous (almost certainty) confidence bands defined in (4.2) include the pointwise confidence intervals. See Remark 1.7 (pp. 233–235) in [13] for more details.*

### Illustration: a simple simulation study

In this paragraph, we present some results from a simulation study. We worked with a sample size $n = 1000$, and considered the case where $\mathbf{X} = (X_1, X_2) \in \mathbb{R}^2$ (i.e. $d = 2$) was such that $X_1 \sim \mathcal{U}(-1, 1)$ and $X_2 \sim \mathcal{U}(-1, 1)$, where $\mathcal{U}(a, b)$ stands for the uniform law on (a,b). Set $m_1(x) = 0.5 \times \cos^2(x)$ and $m_2(x) = 0.5 \times \sin^2(x)$. We selected $\psi = \mathbb{1}_{\{. \leq 0.9\}}$, and considered the model $\mathbb{E}[\psi(Y) | X_1 = x_1, X_2 = x_2] = m_1(x_1) + m_2(x_2)$. Under this model, the variable $Y$ was simulated as follows. For each integer $1 \leq i \leq n$, let $p_i = m_1(x_{1,i}) + m_2(x_{2,i})$ where $x_{j,i}$ is the $i$-th observed value of the variable $X_j$, $j = 1, 2$. Note that $0 < p_i < 1$ for every $1 \leq i \leq n$. Each $Y_i$ was then generated as one $\mathcal{U}(0.9 - p_i, 1 + 0.9 - p_i)$ variable. Following this procedure ensured that $\mathbb{P}(Y_i \leq 0.9 | X_i = x_i) = p_i = m_1(x_{1,i}) + m_2(x_{2,i})$. Regarding the censoring variable, we generated an i.i.d. sample $C_1, \ldots, C_n$ such that $C_i \sim \mathcal{U}(0, 1)$. This choice yielded, *a posteriori*, $\mathbb{P}(\delta = 1) \simeq 0.2$. We used Epanechnikov kernels (for $K$, $K_1$ and $K_2$) and selected $q_1 = q_2 = 0.5 \times \mathbb{1}_{[-1,1]}$ (in such a way that the additive component to estimate were $\eta_{\psi,j} = m_j - 0.25$, $j = 1, 2$). As for the bandwidth choice, we opted *a priori* for $h_{1000} = h_{1,1000} = h_{2,1000} = 0.1$. Results are presented in Figure 1. The confidence bands appear to be adequate, in the sense that they contain the true value of the additive component for "almost" every $x \in [-1, 1]$. The fact that the true function does not belong to our bands for some points was expected: it is due to the $\varepsilon$ term in (4.1). In other respect, the boundary effect pertaining to kernel estimators is perceptible on the plots of Figure 1. In view of the assumption $(C.4)$, we shall however recall that our theorems do not allow to build confidence bands on the entire $[-1, 1]$, and the plots should only be considered on, typically, $[-0.9, 0.9]$.

## 5. Proofs

### 5.1. Proof of Theorem 3.1

Only the proof for the first component is provided. The proof for the $d - 1$ remaining components follows from similar arguments and is therefore omitted.

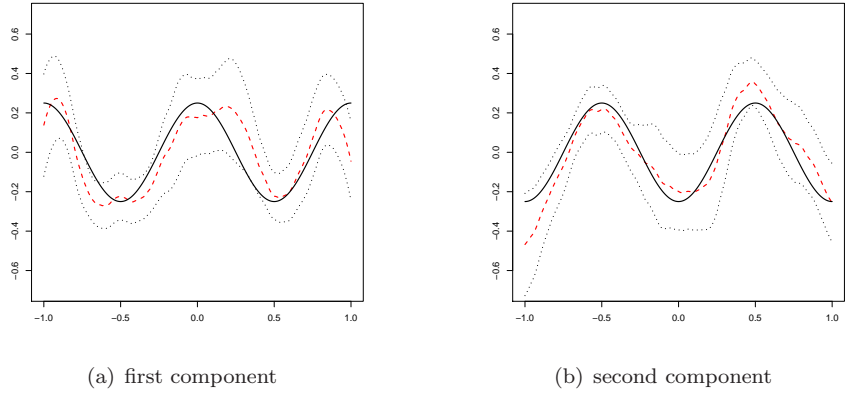

(a) first component                    (b) second component

Fig 1. *Results of the simulation study: true additive components (solid line), their estimates (red dashed line), and the associated confidence bands (dotted line).*

As already mentioned, we first consider the case where both the survival function $G$ of $C$ and the density function $f$ of $\mathbf{X}$ are known.

### 5.1.1. *The case where both $f$ and $G$ are known*

Recall the definitions (2.2) and (2.6) and let $\widetilde{\widehat{m}}_{\psi,n}$ [resp. $\widehat{\widehat{\eta}}_{\psi,1}$] be the version of $\widetilde{m}^\star_{\psi,n}(\mathbf{x})$ [resp. $\widehat{\eta}^\star_{\psi,1}$] in the case where both $G$ and $f$ are known. Namely we have

$$\widetilde{\widehat{m}}_{\psi,n}(\mathbf{x}) \;=\; \frac{1}{n}\sum_{i=1}^{n}\left\{\frac{\delta_i\psi(Z_i)}{G(Z_i)f(\mathbf{X}_i)}\prod_{\ell=1}^{d}\frac{1}{h_{\ell,n}}K_\ell\left(\frac{x_\ell - X_{i,\ell}}{h_{\ell,n}}\right)\right\}, \qquad (5.1)$$

$$\widehat{\widehat{\eta}}_{\psi,1}(x_1) \;=\; \int_{\mathbf{R}^{d-1}}\widetilde{\widehat{m}}_{\psi,n}(\mathbf{x})q_{-1}(\mathbf{x}_{-1})d\mathbf{x}_{-1} - \int_{\mathbf{R}^d}\widetilde{\widehat{m}}_{\psi,n}(\mathbf{x})q(\mathbf{x})d\mathbf{x}. \quad (5.2)$$

In this paragraph, we intend to prove the following result, which is the version of Theorem 3.1 in the case where both $f$ and $G$ are known.

**Proposition 5.1.** *Under the hypotheses of Theorem 3.1, we have,*

$$\sqrt{\frac{nh_{1,n}}{2|\log h_{1,n}|}}\sup_{\psi\in\mathcal{F}}\sup_{x_1\in\mathcal{C}_1}\pm\{\widehat{\widehat{\eta}}_{\psi,1}(x_1) - \eta_{\psi,1}(x_1)\}\xrightarrow{a.s.}\sigma_1,\quad \text{as } n\to\infty, \quad (5.3)$$

*where $\sigma_1$ is as in (2.10).*

In a first step, we will establish Lemma 5.1 below.

**Lemma 5.1.** *Under the assumptions of Theorem 3.1, we have*

$$\sqrt{\frac{nh_{1,n}}{2|\log h_{1,n}|}} \sup_{\psi \in \mathcal{F}} \sup_{x_1 \in \mathcal{C}_1} \pm \{\widehat{\widehat{\eta}}_{\psi,1}(x_1) - \mathbb{E}\widehat{\widehat{\eta}}_{\psi,1}(x_1)\} \xrightarrow{a.s.} \sigma_1, \quad as \ n \to \infty,$$

*where $\sigma_1$ is as in (2.10).*

Let $\psi \in \mathcal{F}$ be a fixed real-valued, measurable and uniformly bounded function defined in $\mathbb{R}$. Following the ideas developed in [15], we will first establish Lemma 5.2 below, which corresponds to Lemma 5.1 in the case where $\mathcal{F}$ is reduced to $\{\psi\}$. Then, we will show how to handle the uniformity over the whole class $\mathcal{F}$ (see Lemma 5.6 and 5.7 in the sequel).

**Lemma 5.2.** *Under the assumptions of Theorem 3.1, we have*

$$\sqrt{\frac{nh_{1,n}}{2|\log h_{1,n}|}} \sup_{x_1 \in \mathcal{C}_1} \pm \{\widehat{\widehat{\eta}}_{\psi,1}(x_1) - \mathbb{E}\widehat{\widehat{\eta}}_{\psi,1}(x_1)\} \xrightarrow{a.s.} \sigma_{\psi,1}, \quad as \ n \to \infty.$$

*where $\sigma_{\psi,1}$ is as in (2.10).*

The proof of Lemma 5.2 is built on recent developments in empirical process theory (see, e.g., [13, 15, 16]). Denote by $\alpha_n$ the multivariate empirical process based upon $(\mathbf{X}_1, Z_1, \delta_1), \ldots, (\mathbf{X}_n, Z_n, \delta_n)$ and indexed by a class $\mathcal{G}$ of measurable functions defined on $\mathbb{R}^{d+2}$. More formally, for $g \in \mathcal{G}$, $\alpha_n(g)$ is defined by

$$\alpha_n(g) = \frac{1}{\sqrt{n}} \sum_{i=1}^{n} \big( g(\mathbf{X}_i, Z_i, \delta_i) - \mathbb{E}g(\mathbf{X}_i, Z_i, \delta_i) \big). \tag{5.4}$$

For $\mathbf{X}_i = (X_{i,1}, \ldots, X_{i,d})$, $1 \le i \le n$, and $x_1 \in \mathcal{C}_1$, set

$$g_{\psi,n}^{x_1}(\mathbf{X}_i, Z_i, \delta_i) = \frac{\delta_i \psi(Z_i)}{G(Z_i)} T_n(\mathbf{X}_i) K_1\left(\frac{x_1 - X_{i,1}}{h_{1,n}}\right), \tag{5.5}$$

$$\text{with} \quad T_n(\mathbf{X}_i) = \frac{1}{f(\mathbf{X}_i)} \int_{\mathbb{R}^{d-1}} \prod_{\ell \neq 1} \frac{1}{h_{\ell,n}} K_\ell\left(\frac{x_\ell - X_{i,\ell}}{h_{\ell,n}}\right) q_\ell(x_\ell) d\mathbf{x}_{-1}. \tag{5.6}$$

From (5.1), (5.2), (5.4) and (5.5), we successively get the two following equalities

$$\sqrt{n}\alpha_n(g_{\psi,n}^{x_1}) = nh_{1,n} \int_{\mathbb{R}^{d-1}} \{\widetilde{\widehat{m}}_{\psi,n}(\mathbf{x}) - \mathbb{E}\widetilde{\widehat{m}}_{\psi,n}(\mathbf{x})\} q_{-1}(x_{-1}) d\mathbf{x}_{-1}, \tag{5.7}$$

$$nh_{1,n}\{\widehat{\widehat{\eta}}_{\psi,1}(x_1) - \mathbb{E}\widehat{\widehat{\eta}}_{\psi,1}(x_1)\} = \sqrt{n}\bigg\{\alpha_n(g_{\psi,n}^{x_1}) - \int_{\mathbb{R}} \alpha_n(g_{\psi,n}^{x_1}) q_1(x_1) dx_1\bigg\}. \tag{5.8}$$

Lemma 5.3 below enables to evaluate the respective order of each of the terms in the right hand side of (5.8). Its proof is postponed until Section 5.2.

**Lemma 5.3.** *Under the conditions of Theorem 3.1, we have, almost surely,*

$$\int_{\mathbb{R}} \alpha_n(g_{\psi,n}^{x_1}) q_1(x_1) dx_1 = o\Big( \sup_{x_1 \in \mathcal{C}_1} \alpha_n(g_{\psi,n}^{x_1}) \Big), \text{ as } n \to \infty. \tag{5.9}$$

In view of (5.8) and (5.9), the asymptotic behavior of the left hand side of (5.8) can be deduced from that of $\alpha_n(g_{\psi,n}^{x_1})$. Then, following once again the ideas of [15] (see also [13]), the proof of Lemma 5.2 will be split into an upper bound part (captured in Lemma 5.4) and a lower bound part (captured in Lemma 5.5).

**Upper bound part**

**Lemma 5.4.** *Recall the definitions (2.10), (5.4) and (5.5). Under the assumptions of Theorem 3.1, we have, for all $\varepsilon > 0$, with probability one,*

$$\limsup_{n \to \infty} \frac{\sup_{x_1 \in \mathcal{C}_1} |\alpha_n(g_{\psi,n_k}^{x_1})|}{\sqrt{2h_{1,n} |\log h_{1,n}|}} \leq (1 + 2\varepsilon)\sigma_{\psi,1}. \tag{5.10}$$

PROOF. We will first examine the behavior of the process $\alpha_n(g_{\psi,n}^{x_1})$ on an appropriately chosen grid of $\mathcal{C}_1$ (partitioning). To do so, we will make use of Bernstein's maximal inequality. Then, we will evaluate the uniform oscillations of our process between the grid points (evaluation of the oscillations). Towards this aim, we will make use of an inequality due to Mason [33], recalled for convenience in Inequality A.1 (see the Appendix).

*Partitioning.* Let $a_1$ and $c_1$ be such that $\mathcal{C}_1 = [a_1, c_1]$, and fix $0 < \delta < 1$. From now on, set, for some $\lambda > 1$, $n_k = [\lambda^k]$, for all $k \geq 1$, and consider the following partitioning of the compact $\mathcal{C}_1$,

$$x_{1,j} = a_1 + j\delta h_{1,n_k}, \quad 0 \leq j \leq J_k := \left[ \frac{c_1 - a_1}{\delta h_{1,n_k}} \right], \tag{5.11}$$

where $u \leq [u] < u + 1$ denotes the integer part of $u$.

Here, we claim that, for all $\varepsilon > 0$, with probability one,

$$\limsup_{k \to \infty} \frac{\max_{n_{k-1} \leq n \leq n_k} \max_{1 \leq j \leq J_k} |\sqrt{n}\alpha_n(g_{\psi,n_k}^{x_{1,j}})|}{\sqrt{2n_k h_{1,n_k} |\log h_{1,n_k}|}} \leq (1 + \varepsilon)\sigma_{\psi,1}. \tag{5.12}$$

For any real valued function $\varphi$ defined on a set $B$, we use the notation $\|\varphi\|_B = \sup_{x \in B} |\varphi(x)|$, and in the particular case where $B = \mathbb{R}^m$, for $m \geq 1$, we will write $\|\varphi\| = \|\varphi\|_B$. Recall that $K_\ell, \ell = 1, \dots, d$, is of bounded variation, and $\psi$ is uniformly bounded. Thus, under the hypothesis (**A**), there exists a constant $0 < \kappa < \infty$ such that, for each $0 \leq j \leq J_k$ and any $x_1 \in \mathcal{C}_1$,

$$\|g_{\psi,n_k}^{x_{1,j}}\| + \|g_{\psi,n}^{x_1}\| \leq \kappa. \tag{5.13}$$

Moreover, by $(C.1)$, and making use of a classical conditioning argument, it follows from (2.8), (5.5) and (5.6) that, for all $0 \leq j \leq J_k$, $k \geq 1$, $1 \leq i \leq n$,

$n_{k-1} < n \leq n_k$,

$$
\begin{aligned}
\mathrm{Var}\big[g_{\psi,n_k}^{x_{1,j}}(\mathbf{X}_i, Z_i, \delta_i)\big] &\leq \mathbb{E}\Big\{\big(g_{\psi,n_k}^{x_{1,j}}(\mathbf{X}_i, Z_i, \delta_i)\big)^2\Big\} \\
&\leq \mathbb{E}\bigg\{\frac{\delta_i \psi^2(Z_i)}{G^2(Z_i)}T_{n_k}^2(\mathbf{X}_i)K_1^2\bigg(\frac{x_{1,j} - X_{1,1}}{h_{1,n_k}}\bigg)\bigg\} \\
&\leq \mathbb{E}\bigg\{\mathbb{E}\bigg(\frac{\delta_i \psi^2(Y_i)}{G^2(Y_i)}\Big|Y_i, \mathbf{X}_i\bigg)T_{n_k}^2(\mathbf{X}_i)K_1^2\bigg(\frac{x_{1,j} - X_{i,1}}{h_{1,n_k}}\bigg)\bigg\} \\
&\leq \mathbb{E}\bigg\{H_\psi(\mathbf{X}_i)T_{n_k}^2(\mathbf{X}_i)K_1^2\bigg(\frac{x_{1,j} - X_{i,1}}{h_{1,n_k}}\bigg)\bigg\} \\
&\leq \int_{\mathbb{R}^d}\frac{H_\psi(\mathbf{u})}{f(\mathbf{u})}\bigg\{\int_{\mathbb{R}^{d-1}}\prod_{\ell\neq 1}\frac{1}{h_{\ell,n_k}}K_\ell\bigg(\frac{x_\ell - u_\ell}{h_{\ell,n_k}}\bigg)q_\ell(x_\ell)d\mathbf{x}_{-1}\bigg\}^2 \\
&\quad \times K_1^2\bigg(\frac{x_{1,j} - u_1}{h_{1,n_k}}\bigg)d\mathbf{u}. \tag{5.14}
\end{aligned}
$$

But, by setting $\mathbf{h}_{-1} = (h_{2,n_k}, \ldots, h_{d,n_k})^T$ and making use of classical changes of variables, it can be derived that, under $(K.1)$ and $(Q.1)$, for a given $0 < \theta < 1$,

$$
\begin{aligned}
&\int_{\mathbb{R}^{d-1}}\prod_{\ell\neq 1}\frac{1}{h_{\ell,k}}K_\ell\bigg(\frac{x_\ell - u_\ell}{h_{\ell,n_k}}\bigg)q_\ell(x_\ell)d\mathbf{x}_{-1} \\
&= \int_{\mathbb{R}^d}\prod_{\ell=1}^d K_\ell(v_\ell)\bigg[q_{-1}(\mathbf{u}_{-1}) \\
&\quad + \sum_{s_2+\cdots+s_d=s}v_2^{s_2}\ldots v_d^{s_d}h_{2,n_k}^{s_2}\ldots h_{d,n_k}^{s_d}\frac{\partial^s q_{-1}}{\partial v_2^{k_1}\ldots \partial v_d^{k_d}}(\mathbf{v}_{-1}\mathbf{h}_{-1}\theta + \mathbf{u}_{-1})\bigg]d\mathbf{v} \\
&= q_{-1}(\mathbf{u}_{-1}) + o(1),
\end{aligned}
$$

in such a way that

$$
\bigg(\int_{\mathbb{R}^{d-1}}\prod_{\ell\neq 1}\frac{1}{h_{\ell,n_k}}K_\ell\bigg(\frac{x_\ell - u_\ell}{h_{\ell,n_k}}\bigg)q_\ell(x_\ell)d\mathbf{x}_{-1}\bigg)^2 = q_{-1}^2(\mathbf{u}_{-1}) + o(1). \tag{5.15}
$$

Recalling the definition (2.10) of $\sigma_{\psi,1}$, it readily follows, from (5.14) and (5.15), that, for all $\varepsilon > 0$ and for $n$ large enough,

$$
\max_{0\leq j\leq J_k}\mathrm{Var}\big(g_{\psi,k}^{x_{1,j}}(\mathbf{X}_i, Z_i, \delta_i)\big) \leq (1+\varepsilon)\sigma_{\psi,1}^2 h_{1,n_k}. \tag{5.16}
$$

In view of (5.4), (5.13) and (5.16), we can apply Bernstein's maximal inequality (see for instance Lemma 2.2 in [14]) to the sequence of random variables,

$$
g_{\psi,n_k}^{x_{1,j}}(\mathbf{X}_i, Z_i, \delta_i) - \mathbb{E}g_{\psi,n_k}^{x_{1,j}}(\mathbf{X}_i, Z_i, \delta_i), \quad i = 1, \ldots, n.
$$

This yields, for $n$ large enough,

$$\mathbb{P}\Big\{ \max_{n_{k-1} \leq n \leq n_k} \max_{1 \leq j \leq J_k} |\alpha_n(g_{\psi,n_k}^{x_{1,j}})| \geq \sigma_{\psi,1}(1+\varepsilon)\sqrt{2h_{1,n_k}|\log h_{1,n_k}|} \Big\}$$

$$\leq 2(J_k+1)\exp\left( \frac{-2\sigma_{\psi,1}^2(1+\varepsilon)h_{1,n_k}|\log h_{1,n_k}|}{2\sigma_{\psi,1}^2 h_{1,n_k} + \frac{2\kappa\sigma_{\psi,1}}{3\sqrt{n_k}}\sqrt{2h_{1,n_k}|\log h_{1,n_k}|}} \right)$$

$$\leq 2(J_k+1)h_{1,n_k}^{1+\varepsilon/2}. \tag{5.17}$$

Keep in mind the definition (5.11) of $J_k$. Since, under $(H.5)$, $\sum_{k \geq 1} h_{1,n}^\varrho < \infty$, for all $\varrho > 0$, the result (5.17) combined with the Borel-Cantelli Lemma naturally implies (5.12).

*Evaluation of the oscillations.* In the sequel, for any class $\mathcal{G}$ of measurable functions, we will denote by $\|\alpha_n\|_{\mathcal{G}} = \sup_{g \in \mathcal{G}} |\alpha_n(g)|$, with $\alpha_n$ as in (5.4).

For future use, first consider a slightly wider class of functions than the one strictly needed in this paragraph. Namely, set

$$\mathcal{G}'_k = \{g_{\psi_1,n_k}^{x_a} - g_{\psi_2,n}^{x_b}, \ n_{k-1} < n \leq n_k, \ x_a, x_b \in \mathcal{C}_1, \ \psi_1, \psi_2 \in \mathcal{F}\}. \tag{5.18}$$

Arguing exactly as in pages 17 and 18 of [15], it can be shown that, for all $k \geq 1$, $\mathcal{G}'_k$ is included in a class $\mathcal{G}'$ of measurable functions, which has a uniform polynomial covering number, i.e., such that for some $C_0 > 0$ and $\mu > 0$, and all $0 < \varepsilon < 1$, $\mathcal{N}(\varepsilon, \mathcal{G}') \leq C_0 \varepsilon^{-\mu}$. Here $\mathcal{N}(\varepsilon, \mathcal{G}') := \sup\{\mathcal{N}(\varepsilon, \mathcal{G}', L_2(\mathbb{P})), \mathbb{P}$ probability measure$\}$ denotes the uniform covering number of the class $\mathcal{G}'$ for $\varepsilon$ and the class of norms $\{L_2(\mathbb{P})\}$, with $\mathbb{P}$ varying in the set of all probability measures on $\mathbb{R}^{d+2}$ (for more details, see, e.g., pp. 83–84 in [41]).

To study the behavior of the process $\alpha_n(g_{\psi,n}^{x_1})$ between the grid points $x_{1,j}$ and $x_{1,j+1}$, with $0 \leq j \leq J_k - 1$, we introduce the following class of functions

$$\mathcal{G}'_{k,j} = \{g_{\psi,n_k}^{x_{1,j}} - g_{\psi,n}^{x_1}, n_{k-1} < n \leq n_k, \ x_{1,j} \leq x_1 \leq x_{1,j+1}\}.$$

Note that, for every $0 \leq j \leq J_k - 1$, we have $\mathcal{G}'_{k,j} \subseteq \mathcal{G}'_k \subseteq \mathcal{G}'$.

Now we claim that, for any $\varepsilon > 0$, there exist almost surely a $\delta_\varepsilon$ and a $\lambda_\varepsilon$ such that,

$$\limsup_{k \to \infty} \max_{0 \leq j \leq J_k - 1} \frac{\max_{n_{k-1} < n \leq n_k} \|n^{1/2}\alpha_n\|_{\mathcal{G}'_{k,j}}}{\sqrt{2n_k h_{1,n_k}|\log h_{1,n_k}|}} \leq \varepsilon\sigma_{\psi,1}, \tag{5.19}$$

whenever (5.11) holds with $0 < \delta \leq \delta_\varepsilon$, $1 < \lambda \leq \lambda_\varepsilon$ and $n_k = [\lambda]^k$.

To establish (5.19), we will make use of Inequality A.1 (see the Appendix). Towards this aim, first note that, since $K_1$ is of bounded variation, we can write $K_1 = K_{1,1} - K_{1,2}$ where $K_{1,1}$ and $K_{1,2}$ are two non-decreasing functions of bounded variation on $\mathbb{R}$. Clearly, $K_{1,1}$ and $K_{1,2}$ are such that $|K_1|_v = |K_{1,1}|_v + |K_{1,2}|_v$, with $|.|_v$ denoting total variation. Then, for all $0 \leq j \leq J_k - 1$ and

$x_{1,j} \le x_1 \le x_{1,j+1}$, it follows that

$$\left| K_1\left(\frac{x_{1,j} - X_{1,1}}{h_{1,n_k}}\right) - K_1\left(\frac{x_1 - X_{1,1}}{h_{1,n}}\right) \right| = \left| \int_{(x_1 - X_{1,1})/h_{1,n}}^{(x_{1,j} - X_{1,1})/h_{1,n_k}} dK_1(y) \right|$$

$$\le \int_{\mathbb{R}} \left| \mathbb{I}\left\{ \frac{x_{1,j} - X_{1,1}}{h_{1,n_k}} > y \right\} - \mathbb{I}\left\{ \frac{x_1 - X_{1,1}}{h_{1,n}} > y \right\} \right| d\big(K_{1,1}(y) + K_{1,2}(y)\big).$$

Since $|K_1(x) - K_1(y)| \le |K_1|_v$ for all $x, y \in \mathbb{R}$, we get

$$\mathbb{E}\left| K_1\left(\frac{x_{1,j} - X_{1,1}}{h_{1,n_k}}\right) - K_1\left(\frac{x_1 - X_{1,1}}{h_{1,n}}\right) \right|^2$$

$$\le |K|_v \int_{\mathbb{R}} \left| \int_{x_1 - yh_{1,n}}^{x_{1,j} - yh_{1,n_k}} f_1(u_1) du_1 \right| d\big(K_{1,1}(y) + K_{1,2}(y)\big)$$

$$\le \|f_1\|_{\mathcal{C}_1^\alpha} |K|_v^2 \left| \frac{h_{1,n} - h_{1,n_k}}{h_{1,n_k}} + \delta \right| h_{1,n_k}. \tag{5.20}$$

Now setting, for $0 \le j \le J_k - 1$,

$$\sigma^2_{\mathcal{G}'_{k,j}} = \sup_{g \in \mathcal{G}'_{k,j}} \mathrm{Var}(g_\psi(\mathbf{X}, Y, \delta)),$$

and making use of the same arguments as those used to derive (5.16), it is readily shown that

$$\sigma^2_{\mathcal{G}'_{k,j}} \le h_{1,n_k} \frac{\|f_1\|_{\mathcal{C}_1^\alpha} |K|_v^2}{\int_{\mathbb{R}} K_1^2} \left| \frac{h_{1,n} - h_{1,n_k}}{h_{1,n_k}} + \delta \right| \sigma^2_{\psi,1}.$$

Set $\tau = 1/[D_1(1 + \sqrt{2/A_2})]$, where $D_1$ and $A_2$ are the constants involved in Inequality A.1. By selecting $\delta > 0$ sufficiently small, and $\lambda > 1$ close enough to 1 to make $\max_{n_{k-1} < n \le n_k} |h_n - h_{n_k}|/h_{n_k}$ as small as desired for large $k$ (using $(H.1\text{-}2)$), we get

$$\sigma^2_{\mathcal{G}'_{k,j}} \le \tau^2 \sigma^2_{\psi,1} \varepsilon^2 h_{1,n_k}. \tag{5.21}$$

Now observe that for all $0 \le j \le J_k - 1$, we have $\|g_\psi\| \le \kappa$ uniformly over $g_\psi \in \mathcal{G}'_{k,j} \subset \mathcal{G}'_k$, where $\kappa$ is as in (5.13). Therefore, applying Inequality A.1 with $\tau$ as in (5.21) and $\rho = \tau\sqrt{2/A_2}$ yields

$$\mathbb{P}\left\{ \max_{0 \le j \le J_k - 1} \frac{\max_{n_{k-1} < n \le n_k} \|n^{1/2}\alpha_n\|_{\mathcal{G}'_{k,j}}}{\sqrt{n_k h_{1,n_k} \log(1/h_{1,n_k})}} \ge \varepsilon \right\} \le 3J_k h_{1,n_k}^2. \tag{5.22}$$

Arguing as before, (5.19) now follows under $(H.5)$ from (5.22) and the definition (5.11) of $J_k$, in combination with the Borel-Cantelli Lemma.

*Conclusion*: The proof of Lemma 5.4 is completed by combining (5.12) and (5.19). $\qquad\square$

**Lower bound part**

**Lemma 5.5.** *Recall the definitions (5.4), (5.5) and (2.10). Under the assumptions of Theorem 3.1, we have, with probability one,*

$$\liminf_{n \to \infty} \sup_{x_1 \in \mathcal{C}_1} \frac{|\alpha_n(g_n^{x_1})|}{\sqrt{2h_{1,n} \log (1/h_{1,n})}} \geq \sigma_{\psi,1}. \qquad (5.23)$$

PROOF. Recall the definition (2.9), and note that, from Scheffe's Lemma, it follows under (**A**) and (*C*.2-3) that the function

$$x_1 \longrightarrow \sqrt{\frac{\phi_{\psi,1}(x_1)}{f_1(x_1)}} \left[ \int_{\mathbb{R}} K_1^2 \right]^{1/2}$$

is continuous on $\mathcal{C}_1$ (see Section $A$.3 in [13] for a complete proof of such continuity results). Then, for any $\epsilon > 0$, we can select a sub-interval $J = [a', c'] \subset \mathcal{C}_1$, such that $\mathbb{P}\{X_1 \in J\} \leq 1/2$ and

$$\inf_{u_1 \in J} \sqrt{\frac{\phi_{\psi,1}(u_1)}{f_1(u_1)}} \left[ \int_{\mathbb{R}} K_1^2 \right]^{1/2} > \sigma_{\psi,1}(1 - \epsilon/2).$$

Now, consider the following partitioning of $J$

$$x_{1,i} = A + 2jh_n, \text{ for } i = 1, \ldots, [(B-A)/2h_{1,n}] - 1 =: k_n.$$

For each $x_{1,i}$, $1 \leq i \leq k_n$, define the function

$$g_i^{(n)}(\mathbf{x}, y, c) = \begin{cases} \dfrac{\psi(y)}{G(y)} T_n(\mathbf{x}) K_1 \left( \dfrac{x_{1,i} - x_1}{h_{1,n}} \right) & \text{if } y \leq c, \\ 0 & \text{if } y > c, \end{cases}$$

where $T_n$ is as in (5.6). Given these notations, the proof of Lemma 5.5 follows from the same lines as those used to establish Proposition 3 in [15]. For the sake of brevity, we omit the details of these book-keeping arguments. $\qquad \square$

From Lemmas 5.3, 5.4 and 5.5, we achieve the proof of Lemma 5.2.

Under the conditions of Theorem 3.1, we readily obtain from Lemma 5.2 that, with probability one, for any finite subclass $\mathcal{G} \subset \mathcal{F}$,

$$\sqrt{\frac{nh_{1,n}}{2|\log h_{1,n}|}} \sup_{\psi \in \mathcal{G}} \sup_{x_1 \in \mathcal{C}_1} \pm \{\widehat{\widehat{\eta}}_{\psi,1}(x_1) - \mathbb{E}\widehat{\widehat{\eta}}_{\psi,1}(x_1)\} \xrightarrow{\text{a.s.}} \sigma_1 \quad \text{as } n \to \infty. \quad (5.24)$$

Therefore, to achieve the proof of Lemma 5.1 we shall show how to extend (5.24) to the entire class $\mathcal{F}$. The following couple of lemmas are directed towards this aim.

**Lemma 5.6.** *Assume the assumptions of Theorem 3.1 hold. For all $\varepsilon > 0$, we can find a finite subclass $\mathcal{G}_\varepsilon \subset \mathcal{F}$, such that, for any $\psi_1 \in \mathcal{F}$, and for $n$ large enough,*

$$\min_{\psi_2 \in \mathcal{G}_\varepsilon} \sup_{x_1 \in \mathcal{C}_1} \mathbb{E}\{[g^{x_1}_{\psi_1,n}(\mathbf{X}, Z, \delta) - g^{x_1}_{\psi_2,n}(\mathbf{X}, Z, \delta)]^2\} \leq \varepsilon h_{1,n},$$

*where, for all $\psi \in \mathcal{F}$, $g^{x_1}_{\psi,n}$ is as in (5.5).*

PROOF. Set $\omega_0 = \omega$ [resp. $\omega_0 = T_F < \infty$] if $(\mathbf{A})(i)$ [resp. $(\mathbf{A})(ii)$] holds. Under $(\mathbf{A})$, it follows from (5.5) and (5.6) that, for $\psi_1, \psi_2 \in \mathcal{F}$ and $x_1 \in \mathcal{C}_1$,

$$\begin{aligned}
&\mathbb{E}\{[g^{x_1}_{\psi_1,n}(\mathbf{X}, Z, \delta) - g^{x_1}_{\psi_2,n}(\mathbf{X}, Z, \delta)]^2\} \\
&= \mathbb{E}\left\{\left[\frac{\mathbb{1}_{\{Y \leq C\}}}{G(Y)} T_n(\mathbf{X}) K_1\left(\frac{x_1 - X_1}{h_{1,n}}\right)(\psi_1(Y) - \psi_2(Y))\right]^2\right\} \\
&\leq \beta h_{1,n} \int_0^{\omega_0} [\psi_1(y) - \psi_2(y)]^2 dy,
\end{aligned}$$

*where* $\beta = G^{-2}(\omega_0)\|K_1\|^2 \sup_{(\mathbf{x},y) \in \mathcal{C}^\alpha \times [0,\omega_0]}\{f_{\mathbf{X},Y}(\mathbf{x}, y)T_n(\mathbf{x})\} < \infty$. *Besides, since $\mathcal{F}$ is a $VC$ subgraph class, it is totally bounded with respect to $d_Q$, where $Q$ is the uniform $(0, \omega_0)$ distribution. Thus, for any $\varepsilon > 0$, we can find a finite class $\mathcal{G}_\varepsilon$ such that*

$$\sup_{\psi_1 \in \mathcal{F}} \min_{\psi_2 \in \mathcal{G}_\varepsilon} \int_0^{\omega_0} [\psi_1(y) - \psi_2(y)]^2 dy \leq \varepsilon/\beta.$$

□

Fix $\varepsilon > 0$ and select $n_0 > 0$ so large that (5.25) holds for all $n \geq n_0$. Further define, for all $\psi_1, \psi_2 \in \mathcal{F}$,

$$d^2(\psi_1, \psi_2) = \sup_{n \geq n_0} h_{1,n}^{-1} \sup_{x_1 \in \mathcal{C}_1} \mathbb{E}\{[g^{x_1}_{\psi_1,n}(\mathbf{X}, Z, \delta) - g^{x_1}_{\psi_2,n}(\mathbf{X}, Z, \delta)]^2\}.$$

Now consider the class of functions

$$\mathcal{G}_n(\varepsilon) = \{g^{x_1}_{\psi_1,n} - g^{x_1}_{\psi_2,n}, d^2(\psi_1, \psi_2) \leq \varepsilon, x_1 \in \mathcal{C}_1\}.$$

**Lemma 5.7.** *Under the assumptions of Theorem 3.1, we have, with probability one,*

$$\limsup_{n \to \infty} \frac{\sup_{d^2(\psi_1,\psi_2) \leq \varepsilon} \sup_{x_1 \in \mathcal{C}_1} \|\alpha_n\|_{\mathcal{G}_n(\varepsilon)}}{\sqrt{2h_{1,n}|\log h_{1,n}|}} \leq A\sqrt{\varepsilon}, \qquad (5.25)$$

*where $A$ is an absolute constant.*

PROOF. The proof of (5.25) is similar to that of (5.19). Set $n_k = 2^k$ and note that,

$$\max_{n_{k-1} < n \leq n_k} \|\alpha_n\|_{\mathcal{G}_n(\varepsilon)} \leq \max_{1 \leq n \leq n_k} \|\alpha_n\|_{\tilde{\mathcal{G}}_k(\varepsilon)},$$

where $\tilde{\mathcal{G}}_k(\varepsilon) = \cup_{n_{k-1}+1}^{n_k} \mathcal{G}_n(\varepsilon)$. It is straightforward that $\sup_{v \in \tilde{\mathcal{G}}_k(\varepsilon)} \|v\| \leq \kappa$, with $\kappa$ as in (5.13). Moreover, keeping in mind the definition (5.18) of $\mathcal{G}'_k$, we have $\mathcal{G}_n(\varepsilon) \subset \mathcal{G}'_k$. Next, observe that, for all large $k$, under $(H.2)$,

$$\sup_{v \in \tilde{\mathcal{G}}_k(\varepsilon)} \text{Var}(v(\mathbf{X}, Z, \delta)) \leq \varepsilon h_{1,n_{k-1}} \leq 2\varepsilon h_{1,n_k}.$$

Arguing as before, Inequality A.1, when applied with $\tau = \sqrt{2\varepsilon}$ and $\rho = \tau\sqrt{1/A_2}$, enables to complete the proof of Lemma 5.7. $\qquad\square$

Recalling (5.8), the proof of Lemma 5.1 is achieved by combining (5.24) with the results of Lemma 5.6 and 5.7. Now, to conclude the proof of Proposition 5.1, it is clearly enough to establish the following result.

**Lemma 5.8.** *Recall the definition (5.2). Under the assumptions of Theorem 3.1, we have,*

$$\sup_{x_1 \in \mathcal{C}_1} \sup_{\psi \in \mathcal{F}} \frac{\sqrt{nh_{1,n}} \left\{ \mathbb{E}\widehat{\eta}_{\psi,1}(x_1) - \eta_{\psi,1}(x_1) \right\}}{\sqrt{|\log h_{1,n}|}} = o(1) \quad as \ \ n \to \infty.$$

PROOF. From (5.1), and arguing as in (2.1), it holds that

$$
\begin{aligned}
\mathbb{E}\widetilde{m}_{\psi,n}(\mathbf{x}) &= \mathbb{E}\left\{ \frac{\delta\psi(Z)}{G(Z)f(\mathbf{X})} \prod_{\ell=1}^{d} \frac{1}{h_{\ell,n}} K_\ell\left(\frac{x_\ell - X_\ell}{h_{\ell,n}}\right) \right\} \\
&= \mathbb{E}\left\{ \frac{\mathbb{E}(\psi(Y)|\ \mathbf{X})}{f(\mathbf{X})} \prod_{\ell=1}^{d} \frac{1}{h_{\ell,n}} K_\ell\left(\frac{x_\ell - X_\ell}{h_{\ell,n}}\right) \right\} \\
&= \int_{\mathbf{R}^d} m_\psi(\mathbf{u}) \prod_{\ell=1}^{d} \frac{1}{h_{\ell,n}} K_\ell\left(\frac{x_\ell - u_\ell}{h_{\ell,n}}\right) d\mathbf{u}.
\end{aligned}
$$

Then, by making use a Taylor development of order $s$ (rendered possible by the assumptions $(K.1)$ and $(C.3)$), we get

$$\sup_{\mathbf{x} \in \mathcal{C}} \sup_{\psi \in \mathcal{F}} |\mathbb{E}\widetilde{m}_{\psi,n}(\mathbf{x}) - m_\psi(\mathbf{x})| = \mathcal{O}\left( \prod_{\ell=1}^{d} h_{\ell,n}^{s_\ell} \right). \tag{5.26}$$

By $(H.3)$, the result of Lemma 5.8 is now a direct consequence of (5.2). $\qquad\square$

### 5.1.2. *Two useful approximation lemmas*

Now, we shall show how to treat the general case (i.e. when neither $f$ nor $G$ is known). Let $\widetilde{m}_{\psi,n}$ [resp. $\widehat{\eta}_{\psi,1}$] be the version of $\widetilde{m}^\star_{\psi,n}(\mathbf{x})$ [resp. $\widehat{\eta}^\star_{\psi,1}$] (see (2.2) [resp. (2.6)]) in the case where $G$ is known and $f$ is unknown. Namely, we have

$$\widetilde{m}_{\psi,n}(\mathbf{x}) = \frac{1}{n}\sum_{i=1}^{n} \left\{ \frac{\delta_i\psi(Z_i)}{G(Z_i)\widehat{f}_n(\mathbf{X}_i)} \prod_{\ell=1}^{d} \frac{1}{h_{\ell,n}} K_\ell\left(\frac{x_\ell - X_{i,\ell}}{h_{\ell,n}}\right) \right\}, \tag{5.27}$$

$$\widehat{\eta}_{\psi,1}(x_1) = \int_{\mathbf{R}^{d-1}} \widetilde{m}_{\psi,n}(\mathbf{x})q_{-1}(\mathbf{x}_{-1})d\mathbf{x}_{-1} - \int_{\mathbf{R}^d} \widetilde{m}_{\psi,n}(\mathbf{x})q(\mathbf{x})d\mathbf{x}. \tag{5.28}$$

**Lemma 5.9.** *Recall the definition (5.2). Under the assumptions of Theorem 3.1, we have, almost surely,*

$$\sup_{x_1 \in \mathcal{C}_1} \sup_{\psi \in \mathcal{F}} \frac{\sqrt{nh_{1,n}}\{\widehat{\eta}_{\psi,1}(x_1) - \widehat{\widehat{\eta}}_{\psi,1}(x_1)\}}{\sqrt{|\log h_{1,n}|}} = o(1) \quad as \ n \to \infty. \tag{5.29}$$

PROOF. Because $q(\mathbf{x}) = \prod_{\ell=1}^d q_\ell(x_\ell)$ and because the functions $q_\ell$, $\ell = 1, \ldots, d$, are bounded, a classical change of variable yields, for $i = 1, \ldots n$,

$$\int_{\mathbb{R}^{d-1}} \left| \prod_{\ell=1}^d \frac{1}{h_{\ell,n}} K_\ell\left(\frac{x_\ell - X_{i,\ell}}{h_{\ell,n}}\right) q(\mathbf{x})d\mathbf{x} \right| \leq M_1,$$

with $0 < M_1 < \infty$. Recall the definition (5.1) and set $\Psi(y,c) = \mathbb{I}_{\{y \leq c\}}\psi(y \wedge c)/G(y \wedge c)$, for all $y, c \in \mathbb{R}$. Then, since, for $\ell = 1, \ldots, d$, $q_\ell$ has a compact support included in $\mathcal{C}_\ell$ and $K_\ell$ is compactly supported, we have under $(H.1)$ and for $n$ large enough,

$$\int_{\mathbb{R}^d} \left| \widetilde{m}_{\psi,n}(\mathbf{x}) - \widetilde{\widetilde{m}}_{\psi,n}(\mathbf{x}) \right| q(\mathbf{x})d\mathbf{x} \leq \frac{M_1}{n} \sum_{i=1}^n \left| \Psi(Y_i, C_i) \right| \sup_{\mathbf{x} \in \mathcal{C}^\alpha} \frac{|\hat{f}_n(\mathbf{x}) - f(\mathbf{x})|}{|f(\mathbf{x})\hat{f}_n(\mathbf{x})|}.$$

Clearly, by $(\mathbf{A})$, $\Psi$ is uniformly bounded. Therefore, the following result (see, e.g., [1])

$$\sup_{\mathbf{x} \in \mathcal{C}^\alpha} |\hat{f}_n(\mathbf{x}) - f(\mathbf{x})| = \mathcal{O}\left(\sqrt{\frac{\log n}{nh_n^d}}\right) \quad \text{a.s., as } n \to \infty,$$

is enough to conclude under $(C.4)$ that, almost surely as $n \to \infty$,

$$\int_{\mathbb{R}^d} \left| \widetilde{m}_{\psi,n}(\mathbf{x}) - \widetilde{\widetilde{m}}_{\psi,n}(\mathbf{x}) \right| q(\mathbf{x})d\mathbf{x} = \mathcal{O}\left(\sqrt{\frac{\log n}{nh_n^d}}\right).$$

Similarly, it can be shown that, almost surely as $n \to \infty$,

$$\sup_{x_1 \in \mathcal{C}_1} \int_{\mathbb{R}^{d-1}} \left| \widetilde{m}_{\psi,n}(x_1, \mathbf{x}_{-1}) - \widetilde{\widetilde{m}}_{\psi,n}(x_1, \mathbf{x}_{-1}) \right| q_{-1}(\mathbf{x}_{-1})d\mathbf{x}_{-1} = \mathcal{O}\left(\sqrt{\frac{\log n}{nh_n^d}}\right).$$

From these two last statements and the definitions (5.2) and (5.28), the proof of Lemma 5.9 is completed under the assumption $(H.4)$. □

**Lemma 5.10.** *Recall the definitions (2.6) and (5.28). Under the assumptions of Theorem 3.1, we have, almost surely as $n \to \infty$,*

$$\sup_{x_1 \in \mathcal{C}_1} \sup_{\psi \in \mathcal{F}} \frac{\sqrt{nh_{1,n}}\{\widehat{\eta}^\star_{\psi,1}(x_1) - \widehat{\eta}_{\psi,1}(x_1)\}}{\sqrt{|\log h_{1,n}|}} = o(1). \tag{5.30}$$

PROOF. First consider the case where $(\mathbf{A})(i)$ holds. Set $b = \inf_{\mathbf{x} \in \mathcal{C}^\alpha} f(\mathbf{x})$. Note that $b > 0$ by $(C.4)$. Then, recalling (2.2) and (5.27) and arguing as we did along the proof of Lemma 5.9, we get

$$\int_{\mathbb{R}^d} \big| \widetilde{m}^\star_{\psi,n}(\mathbf{x}) - \widetilde{m}_{\psi,n}(\mathbf{x}) \big| q(\mathbf{x}) d\mathbf{x} \le \frac{M_1}{b + o(1)} \sup_{0 \le y \le \omega} \left\{ \frac{|\psi(y)| \cdot |G^\star_n(y) - G(y)|}{|G^\star_n(y) G(y)|} \right\}.$$

Since $\omega < T_H$, the iterated law of the logarithm of [18] ensures that

$$\sup_{y \le \omega} |G^\star_n(y) - G(y)| = \mathcal{O}((\log \log n / n)^{1/2})$$

almost surely as $n \to \infty$. Therefore, it follows under $(C.2\text{-}3\text{-}4)$ that, almost surely as $n \to \infty$,

$$\int_{\mathbb{R}^d} \big| \widetilde{m}^\star_{\psi,n}(\mathbf{x}) - \widetilde{m}_{\psi,n}(\mathbf{x}) \big| q(\mathbf{x}) d\mathbf{x} = \mathcal{O}\left( \sqrt{\frac{\log \log n}{n}} \right). \tag{5.31}$$

In the same spirit it can be shown that, almost surely as $n \to \infty$,

$$\sup_{x_1 \in \mathcal{C}_1} \int_{\mathbb{R}^{d-1}} \big| \widetilde{m}^\star_{\psi,n}(\mathbf{x}) - \widetilde{m}_{\psi,n}(\mathbf{x}) \big| q_{-1}(\mathbf{x}_{-1}) d\mathbf{x}_{-1} = \mathcal{O}\left( \sqrt{\frac{\log \log n}{n}} \right), \tag{5.32}$$

which, by $(H.5)$, completes the proof of Lemma 5.10 under $\mathbf{A}(i)$. In the case where $(\mathbf{A})(ii)$ holds, the proof follows from the same lines as above, making use of either the iterated law of the logarithm of [20] (if $(\mathbf{A})(ii)$ holds with $p = 1/2$) or Theorem 2.1 of [6] (if $(\mathbf{A})(ii)$ holds with $0 < p < 1/2$) instead of the iterated law of the logarithm of [18]. The details are omitted. □

By combining Lemmas 5.9 and 5.10 with Proposition 5.1, we conclude the proof of Theorem 3.1.

### 5.2. *Proof of Lemma 5.3*

Set

$$\Psi(Y_i, C_i) = \mathbb{1}_{\{Y_i \le C_i\}} \psi(Y_i \wedge C_i) / G(Y_i \wedge C_i) = \delta_i \psi(Z_i) / G(Z_i)$$

$$\tilde{\Psi}_n(Y_i, C_i) = \Psi(Y_i, C_i) \int_{\mathbb{R}^{d-1}} \prod_{\ell=2}^d \frac{1}{h_{\ell,n}} K_\ell \left( \frac{x_\ell - X_{i,\ell}}{h_{\ell,n}} \right) \frac{q_{-1}(\mathbf{x}_{-1})}{f(X_{i,-1} | X_{i,1})} d\mathbf{x}_{-1},$$

$$\tilde{g}(x_1) = \mathbb{E}\big( \tilde{\Psi}_n(Y_i, C_i) \big| X_{i,1} = x_1 \big), \tag{5.33}$$

$$\text{and} \quad \beta_1(x_1) = \frac{1}{n h_{1,n}} \sum_{i=1}^n \frac{\tilde{\Psi}_n(Y_i, C_i)}{f_1(X_{i,1})} K_1 \left( \frac{x_1 - X_{i,1}}{h_{1,n}} \right).$$

It follows that,

$$
\begin{aligned}
\mathrm{Var}(\beta_1(x_1)) \quad &= \quad \frac{1}{n^2}\sum_{i=1}^{n}\mathbb{E}\left\{\frac{\tilde{\Psi}_n(Y_i,C_i)}{f_1(X_{i,1})h_{1,n}}K_1\left(\frac{x_1-X_{i,1}}{h_{1,n}}\right)\right\}^2 \\
&\qquad -\frac{1}{n^2}\sum_{i=1}^{n}\left\{\mathbb{E}\left(\frac{\tilde{\Psi}_n(Y_i,C_i)}{f_1(X_{i,1})h_{1,n}}K_1\left(\frac{x_1-X_{i,1}}{h_{1,n}}\right)\right)\right\}^2 \\
&=: \quad \frac{1}{n}\left[\frac{1}{h_{1,n}}\Phi_{1,n}(x_1)-[\Gamma_{1,n}(x_1)]^2\right].
\end{aligned}
$$

Observing that $\tilde{g}$ is uniformly bounded under the assumptions of Theorem 3.1 and making use of some conditioning arguments, it can be shown that

$$
\frac{1}{n}\,\Gamma_{1,n}^2(x_1)\to 0 \quad\text{as}\ \ n\to\infty. \tag{5.34}
$$

Moreover, by using the change of variable $v_1h_{1,n}=x_1-u_1$, we obtain

$$
\begin{aligned}
\Phi_{1,n}(x_1) \quad &= \quad \int_{\mathbb{R}}\frac{K_1^2(v_1)}{f_1(x_1-h_{1,n}v_1)}\,\mathbb{E}(\tilde{\Psi}_n^2(Y_i,C_i)|\ X_{i,1}=x_1-h_{1,n}v_1)dv_1 \\
&= \quad \int_{\mathbb{R}}K_1^2(v_1)\left(\frac{\mathbb{E}(\tilde{\Psi}_n^2(Y_i,C_i)\mid X_{i,1}=x_1-h_1v_1)}{f_1(x_1-h_{1,n}v_1)}-\frac{\phi_{\psi,1}(x_1)}{f_1(x_1)}\right)dv_1 \\
&\qquad +\frac{\phi_{\psi,1}(x_1)}{f_1(x_1)}\int_{\mathbb{R}}K_1^2(v_1)dv_1.
\end{aligned}
$$

But, recalling the definition (2.9) of the function $\phi_{\psi,1}$, the quantity

$$
\left|\frac{\mathbb{E}(\tilde{\Psi}(Y_i,C_i)|\ X_{i,1}=x_1-h_{1,n}v_1)}{f_1(x_1-h_{1,n}v_1)}-\frac{\phi_{\psi,1}(x_1)}{f_1(x_1)}\right|
$$

is clearly bounded under the assumptions $(C.4),(C.6),(K.1)$ and $(Q.1)$. Therefore, Lebesgue's dominated convergence Theorem enables us to conclude that

$$
\Phi_{1,n}(x_1)\to\frac{\phi_{\psi,1}(x_1)}{f_1(x_1)}\int_{\mathbb{R}}K_1^2(v_1)dv_1. \tag{5.35}
$$

Combining (5.34) and (5.35) we obtain, for all $x_1\in\mathcal{C}_1$,

$$
\mathrm{Var}(\beta_1(x_1))=\mathbb{E}\left\{\beta_1(x_1)-\mathbb{E}\beta_1(x_1)\right\}^2=\mathcal{O}(n^{-2s/(2s+1)}).
$$

Then,

$$
\begin{aligned}
\int_{\mathcal{C}_1}\mathrm{Var}(\beta_1(x_1))dx_1 \quad &= \quad \int_{\mathcal{C}_1}\mathbb{E}\left\{\beta_1(x_1)-\mathbb{E}\beta_1(x_1)\right\}^2dx_1 \\
&= \quad \mathbb{E}\left(\int_{\mathcal{C}_1}\left\{\beta_1(x_1)-\mathbb{E}\beta_1(x_1)\right\}^2dx_1\right) \\
&= \quad \mathcal{O}(n^{-2k/(2k+1)}). 
\end{aligned} \tag{5.36}
$$

Recall the definitions (5.1), (5.4), (5.5), (5.6) and (5.33). From (5.7), and using the Cauchy-Schwartz inequality, we obtain,

$$\left| \int_{\mathbb{R}} \frac{\alpha_n(g_{\psi,n}^{x_1})}{\sqrt{2h_{1,n}|\log h_{1,n}|}} \, q_1(x_1)dx_1 \right|$$

$$= \sqrt{\frac{nh_{1,n}}{2|\log(h_{1,n})|}} \left| \int_{\mathbb{R}^d} \{\widetilde{\widetilde{m}}_{\psi,n}(\mathbf{x}) - \mathbb{E}\widetilde{\widetilde{m}}_{\psi,n}(\mathbf{x})\} q(\mathbf{x})d\mathbf{x} \right|$$

$$= \sqrt{\frac{nh_{1,n}}{2|\log(h_{1,n})|}} \left| \int_{\mathcal{C}_1} \{\beta_1(x_1) - \mathbb{E}\beta_1(x_1)\} q_1(x_1)dx_1 \right|$$

$$\leq \sqrt{\frac{nh_{1,n}}{2|\log(h_{1,n})|}} \sqrt{\int_{\mathcal{C}_1} \{\beta_1(x_1) - \mathbb{E}\beta_1(x_1)\}^2 dx_1 \int_{\mathcal{C}_1} q_1^2(x_1)dx_1}. \quad (5.37)$$

From (5.36) and (5.37), it follows that, almost surely as $n \to \infty$,

$$\left| \int_{\mathbb{R}} \frac{\alpha_n(g_{\psi,n}^{x_1})}{\sqrt{2h_{1,n}|\log h_{1,n}|}} \, q_1(x_1)dx_1 \right| = \mathcal{O}\left( \sqrt{\frac{n^{1/(2s+1)}h_{1,n}}{|\log h_{1,n}|}} \right). \quad (5.38)$$

But, from (5.23), we have

$$\sup_{x_1 \in \mathcal{C}_1} \frac{\alpha_n(g_{\psi,n}^{x_1})}{\sqrt{2h_{1,n}|\log h_{1,n}|}} \geq \sigma_{\psi,1}. \quad (5.39)$$

The proof of Lemma 5.3 is readily achieved by combining (5.38) and (5.39) with the condition $(H.3)$. $\square$

### 5.3. Proof of Theorem 3.2

Recall the definitions (2.6) and (2.7) and observe that,

$$\left| \sqrt{\frac{nh_{1,n}}{2|\log h_{1,n}|}} \sup_{\psi \in \mathcal{F}} \sup_{\mathbf{x} \in \mathcal{C}} \pm \{\widetilde{m}_{\psi,add}^{\star}(\mathbf{x}) - m_\psi(\mathbf{x})\} - \sum_{\ell=1}^{d} \sigma_\ell \right|,$$

$$\leq \sum_{\ell=1}^{d} \left| \sqrt{\frac{nh_{1,n}}{2|\log h_{1,n}|}} \sup_{\psi \in \mathcal{F}} \sup_{x_\ell \in \mathcal{C}_\ell} \pm \{\widehat{\eta}_{\psi,\ell}^{\star}(x_\ell) - \eta_{\psi,\ell}^{\star}(x_\ell)\} - \sigma_\ell \right|$$

$$+ \sqrt{\frac{nh_{1,n}}{2|\log h_{1,n}|}} \left| \int_{\mathbb{R}^d} \{\widetilde{m}_{\psi,n}^{\star}(\mathbf{x}) - m_\psi(\mathbf{x})\} q(\mathbf{x})d\mathbf{x} \right|.$$

Under the assumption $(H.3)$ and $(H.4)$, by proceeding as we did along the proof of Lemma 5.3 (see also (5.26) and (5.37)), we get

$$\sqrt{\frac{nh_{1,n}}{2|\log h_{1,n}|}} \int_{\mathbb{R}^d} \{\widetilde{m}_{\psi,n}^{\star}(\mathbf{x}) - m(\mathbf{x})\} q(\mathbf{x})d\mathbf{x} = o(1) \quad \text{a.s.}. \quad (5.40)$$

By combining Theorem 3.1 and the statement (5.40), we complete the proof of Theorem 3.2.

# Appendix

In this section, we present an inequality which was of particular interest for our task. It is due to Mason [33], who derived it from an inequality obtained by Talagrand [39].

INEQUALITY A.1. Let $Z, Z_1, \ldots, Z_n$ be i.i.d. random variables, $n \geq 1$. Denote by $\mathcal{F}$ a class of functions such that

$$\sup_{f \in \mathcal{F}} \mathrm{Var}\left(f(Z)\right) \leq \tau^2 h, \quad \text{with } \tau, h > 0.$$

Assume there exist constants $M$, $C$ and $\nu > 0$, fulfilling, for all $0 < \varepsilon < 1$,

$$\mathcal{N}(\varepsilon, \mathcal{F}) \leq C \varepsilon^{-\nu} \quad \text{and} \quad \sup_{f \in \mathcal{F}, z \in \mathbb{R}^d} |f(z)| \leq M.$$

Choose any $\rho > 0$. Then there exist a universal constant $A_2 > 0$, and a constant $D_1 = D_1(\nu) > 0$, depending only on $\nu$, such that, if $h > 0$ satisfies

$$K_1 := \max\left\{\frac{4M\sqrt{\nu+1}}{\tau}, \frac{M\rho}{\tau^2}\right\} \quad \leq \quad \sqrt{\frac{nh}{|\log h|}},$$

$$K_2 := \min\left\{\frac{1}{\tau^2 M}, \tau^2\right\} \quad \geq \quad h,$$

then we have, with $T_n(g) = \sum_{j=1}^{n} \left\{g(Z_j) - \mathbb{E}(g(Z))\right\}$ for $g \in \mathcal{F}$,

$$\mathbb{P}\left(\sup_{1 \leq m \leq n} \|T_m(\cdot)\|_{\mathcal{F}} \geq (\tau + \rho) D_1 \sqrt{nh|\log h|}\right) \leq 2 \exp\left(-\frac{A_2 \rho^2}{\tau^2}|\log h|\right).$$

# References

[1] P. Ango-Nze and R. Rios. Density estimation in $L^\infty$ norm for mixing processes. *J. Stat. Plann. Inference*, 83(1):75–90, 2000. MR1741444

[2] R. Beran. Nonparametric regression with randomly censored data. In *Technical report*. Univ. California Press, Berkeley, 1981.

[3] E. Brunel and F. Comte. Adaptive nonparametric regression estimation in presence of right censoring. *Math. Methods Stat.*, 15(3):233–255, 2006. MR2278288

[4] J. Buckley and I. James. Linear regression with censored data. *Biometrika*, 66:429–464, 1979.

[5] A. Carbonez, L. Györfi, and E.C. van der Meulen. Partitioning-estimates of a regression function under random censoring. *Statist. Decisions*, 13(1):21–37, 1995. MR1334076

[6] K. Chen and S.H. Lo. On the rate of uniform convergence of the Product-Limit estimator: strong and weak laws. *Ann. Statist.*, 25(3):1050–1087, 1997. MR1447741

[7] M. Csörgő and P. Révész. *Strong Approximation in Probability and Statistics*. Acedemic Press, New York, 1981. MR0666546

[8] D.M. Dabrowska. Nonparametric regression with censored covariates. *J. Multivariate Anal.*, 54(2):253–283, 1995. MR1345539

[9] M. Debbarh. Some uniform limit results in additive regression model. *To appear in Communication in Statistics - Theory and Methods*, 2008.

[10] M. Debbarh and V. Viallon. Mean square convergence for an estimator of the additive regression function under random censorship. *C. R. Acad. Sci. Paris Ser. I*, 344(3):205–210, 2007. MR2292289

[11] M. Debbarh and V. Viallon. Asymptotic normality of the additive regression model components under random censorship. *Preprint. Available at http://arxiv.org/abs/math/0612507*, 2008.

[12] P. Deheuvels and G. Derzko. Uniform consistency for conditinal lifetime distribution estimators under random censorship. In J.L. Auget, N. Balakrishnan, M. Mesbah, and G. Molenberghs, editors, *Advances in Statistical Methods in the Health Sciences: Applications to Cancer and AIDS Studies, Genome Sequence Analysis and Survival Analysis*. Birkhäuser, Boston, 2007.

[13] P. Deheuvels and D.M. Mason. General asymptotic confidence bands based on kernel-type function estimators. *Stat. Inference Stoch. Process.*, 7:225–277, 2004. MR2111291

[14] U. Einmahl and D.M. Mason. Some universal results on the behavior of the increments of partial sums. *Ann. Probab.*, 24:1388–1407, 1996. MR1411499

[15] U. Einmahl and D.M. Mason. An empirical process approach to the uniform consistency of kernel type estimators. *J. Theor. Probab.*, 13:1–13, 2000. MR1744994

[16] U. Einmahl and D.M. Mason. Uniform in bandwidth consistency of kernel-type function estimators. *Ann. Statist.*, 33(3):1380–1403, 2005. MR2195639

[17] J. Fan and I. Gijbels. Censored regression: Local linear approximations and their applications. *J. Am. Stat. Assoc.*, 89(426):560–570, 1994. MR1294083

[18] A. Földes and L. Rejtő. A LIL type result for the product-limit estimator. *Z. Wahrsch. Verw. Gebiete*, 56:75–86, 1981. MR0612161

[19] S. Gross and T.L. Lai. Nonparametric estimation and regression analysis with left-truncated and right-censored data. *J. Am. Stat. Assoc.*, 91(426):1166–1180, 1996. MR1424616

[20] M. Gu and T.L. Lai. Functional laws of the iterated logarithm for the product-limit estimator of a distribution function under random censorship or truncation. *Ann. Probab.*, 18:160–189, 1990. MR1043942

[21] T.J. Hastie and R.J. Tibshirani. *Generalized additive models*. Chapman and Hall, 1990. MR1082147

[22] N.W. Hengartnern and S. Sperlich. Rate optimal estimation with the integration method in the presence of many covariates. *J. Multivar. Anal.*, 95:246–272, 2005. MR2170397

[23] J. Horowitz, J. Klemelä, and E. Mammen. Optimal estimation in additive regression models. *Bernoulli*, 12(2):271–298, 2006. MR2218556

[24] M.C. Jones, S.J. Davies, and B.U. Park. Versions of kernel-type regression

estimators. *J. Am. Stat. Assoc.*, 89:825–832, 1994. MR1294728

[25] E.L. Kaplan and P. Meier. Nonparametric estimation from incomplete observations. *J. Am. Stat. Assoc.*, 53:457–481, 1958. MR0093867

[26] M. Kohler, S. Kul, and K. Máthé. Least squares estimates for censored regression. *Preprint. Available at http://www.mathematik.uni-stuttgart.de/mathA/lst3/kohler/hfm-pub-en.html*, 2006.

[27] M. Kohler, K. Máthé, and M. Pintér. Prediction from randomly right censored data. *J. Multivar. Anal.*, 80(1):73–100, 2002. MR1889835

[28] H. Koul, V. Susarla, and J. Van Ryzin. Regression analysis with randomly right-censored data. *Ann. Statist.*, 9:1276–1288, 1981. MR0630110

[29] T.L. Lai, Z. Ying, and Z. Zheng. Asymptotic normality of a class of adaptive statistics with applications to synthetic data methods for censored data. *J. Multivar. Anal.*, 52:259–279, 1995. MR1323333

[30] S. Leurgans. Linear models, random censoring and synthetic data. *Biometrika*, 74:301–309, 1987. MR0903130

[31] O. Linton and J.P. Nielsen. A kernel method of estimating structured nonparametric regression based on marginal integration. *Biometrika*, 82(1):93–100, 1995. MR1332841

[32] E. Mammen and B.U. Park. A simple smooth backfitting method for additive models. *Ann. Statist.*, 34(5):2252–2271, 2006. MR2291499

[33] D.M. Mason. A uniform functional law of the iterated logarithm for the local empirical process. *Ann. Probab.*, 32(2):1391–1418, 2004. MR2060302

[34] W.K. Newey. Kernel estimation of partial means and a general variance estimator. *Econometric Theory*, 10(2):233–253, 1994. MR1293201

[35] D. Rubin and M. van der Laan. A doubly robust censoring unbiased transformation. *Int. J. Biostat.*, 3(1):Article 4, 2007. MR2306842

[36] C.J. Stone. Additive regression and other nonparametric models. *Ann. Statist.*, 13:689–705, 1985. MR0790566

[37] W. Stute. The law of the iterated logarithm for kernel density estimators. *Ann. Prob.*, 10:414–422, 1982. MR0647513

[38] W. Stute. Nonlinear censored regression. *Statistica Sinica*, 9:1089–1102, 1999. MR1744826

[39] M. Talagrand. Sharper bounds for Gaussian and empirical processes. *Ann. Probab.*, 22:28–76, 1994. MR1258865

[40] D. Tjøstheim and B.H. Auestad. Nonparametric identification of nonlinear time series: projections. *J. Am. Statist. Assoc.*, 89(428):1398–1409, 1994. MR1310230

[41] A.W. van der Vaart and J.A. Wellner. *Weak convergence and empirical processes.* Springer, New York, 1996. MR1385671